\documentclass[preprint,12pt]{elsarticle}

\usepackage{amsmath, amsthm, amssymb}
\usepackage{graphicx}
\usepackage{epstopdf}
\usepackage{comment}
\usepackage{caption}
\usepackage{subcaption}

\newtheorem{proposition}{Proposition}

\renewcommand{\vec}[1]{\mathbf{#1}}
\newcommand{\cvec}[1]{\underline{#1}}
\newcommand{\dx}{\Delta x}
\newcommand{\diag}[1]{\mathrm{diag}(#1)}
\newcommand{\ordo}[1]{\mathcal{O}(\dx^{#1})}

\begin{document}

\begin{frontmatter}
\title{A superconvergent stencil-adaptive SBP-SAT finite difference scheme}
\author[LU]{Viktor Linders\corref{cor1}}
\ead{viktor.linders@math.lu.se}
\cortext[cor1]{Corresponding author}
\author[NASA]{Mark Carpenter}
\ead{mark.h.carpenter@nasa.gov}
\author[LiU,Jo]{Jan Nordstr\"{o}m}
\ead{jan.nordstrom@liu.se}
\address[LU]{Center for Mathematical Sciences,
Lund University,
Lund,
Sweden.}
\address[NASA]{Computational AeroSciences Branch,
NASA Langley Research Center,
Hampton, VA,
United States.}
\address[LiU]{Department of Mathematics,
Computational Mathematics,
Link\"{o}ping University,
Link\"{o}ping,
Sweden}
\address[Jo]{Department of Mathematics and Applied Mathematics, 
University of Johannesburg, 
South Africa}

\begin{abstract}
A stencil-adaptive SBP-SAT finite difference scheme is shown to display superconvergent behavior. Applied to the linear advection equation, it has a convergence rate $\ordo{4}$ in contrast to a conventional scheme, which converges at a rate $\ordo{3}$.
\end{abstract}


\end{frontmatter}


\section{Introduction} \label{sec:introduction}

This note concerns finite difference methods where the stencil coefficients are adapted in time. In \cite{linders2020accurate}, a systematic way of constructing adaptive \emph{central} difference stencils of bandwidth $2n+1$ was presented and applied to first order hyperbolic problems with periodic boundary conditions. It was shown that such stencils were more accurate than conventional finite difference schemes based on Taylor expansions, irrespective of the grid resolution. Further, the advantages of Dispersion Relation Preserving (DRP) schemes\footnote{DRP schemes are finite difference methods optimized in Fourier space; see e.g. \cite{linders2015uniformly}.} were recovered for under-resolved problems. This was achieved without loss of the high convergence rate of conventional stencils, which tends to plague DRP schemes. These adaptive stencils are \emph{superconvergent} in the sense that they offer $\ordo{2n}$ convergence despite having considerably lower formal order of accuracy. Here, $\dx$ represents a characteristic grid spacing.

The observations presented herein are a result of attempting to extend the work in \cite{linders2020accurate} to Summation-By-Parts (SBP) operators. Adjoined with Simultaneous Approximation Terms (SATs), SBP methods offer a systematic approach for stable and accurate discretization of initial boundary value problems \cite{fernandez2014review,svard2014review}. Connections between SBP and DRP have been made in \cite{linders2017summation}.

SBP operators consist of a pair of matrices ($P$,$Q$), such that $D = P^{-1}Q$ approximates the first derivative. By definition, they satisfy the properties
\[
P = P^\top > 0 \quad \text{and} \quad Q + Q^\top = \diag{-1,0,\dots,0,1}.
\]
Operators that utilize a diagonal $P$ enjoy advantageous stability properties when used for spatial as well as temporal discretization \cite{nordstrom2013summation,linders2020properties}, essentially because $P$ then commutes with other diagonal matrices. However, finite difference SBP operators with diagonal $P$ suffer from reduced order of accuracy near boundaries \cite{kreiss1974finite,strand1994summation,linders2018order}. If a central difference stencil of order $2p$ is used in the interior, the boundary stencil is limited to order $p$. This typically leads to a convergence rate of order $p+1$ when SBP-SAT is used to discretize first order energy stable hyperbolic IBVPs \cite{svard2019convergence}.

It is possible to increase the boundary accuracy of SBP operators to $2p-1$ by utilizing so-called block norms, meaning that $P$ has non-diagonal blocks near the boundaries. In case a stable scheme can be constructed, this increases the convergence rate to $2p$. However, it is in general not straightforward to obtain energy estimates with such methods. In this note, we find a different use for block-norm SBP operators as 'target stencils' in an optimization problem.

While our ambition has been to develop a general purpose adaptive SBP-SAT method, we disclaim right away that this goal is reached within this note. Instead, we present an example of an adaptive SBP method based on the conventional SBP(4,2) operator\footnote{SBP(4,2) uses a five-point, fourth order accurate central difference stencil in the interior and four rows of second order accurate, non-symmetric stencils near the boundaries. Its convergence rate for first order energy stable hyperbolic problems is $\ordo{3}$ \cite{svard2019convergence}.} and discuss the various difficulties one faces in extending the method to arbitrary accuracy. A remarkable observation is that for the test problem considered, the observed convergence rate of the adaptive SBP method is $\ordo{4}$, which is one order higher than conventional SBP theory allows. The adaptive method is thus superconvergent.

The remainder of this note is organized as follows: In Section \ref{sec:stability} stability requirements on the adaptive scheme are presented. An outline of the adaptive procedure is given in Section \ref{sec:optimization}. Numerical demonstrations of superconvergence are given in Section \ref{sec:experiments}. Comments on generalizations are made in Section \ref{sec:generalization} followed by conclusions in Section \ref{sec:conclusions}.

 
 \section{Remarks on stability} \label{sec:stability}

We will shortly describe how to construct the adaptive SBP operator and apply it to the linear advection equation with periodic boundary conditions:
\begin{equation} \label{eq:advection}
    \begin{alignedat}{5}
    u_t + u_x &= 0, &&x \in \Omega, \quad &&t \in (0,T], \\
    u(x,0) &= u_0(x), \quad && x \in \Omega. &&
    \end{alignedat}
\end{equation}
Before doing so, we briefly discuss the SBP-SAT discretization and the stability requirements on the method.

We divide $\Omega$ into a set of $K$ non-overlapping blocks/elements, each subdivided into $N+1$ uniformly spaced grid points, and discretize with an SBP operator $D_k = P_k^{-1}Q_k$ on the $k$th block:
\[
\vec{u}^{(k)}_t + D_k \vec{u}^{(k)} = -\frac{1+\theta}{2} P_k^{-1} (u_0^{(k)} - u_N^{(k-1)}) \vec{e}_0 + \frac{1-\theta}{2} P_k^{-1} (u_N^{(k)} - u_0^{(k+1)}) \vec{e}_N.
\]
Here, $u_i^{(j)}$ is to be understood as the $i$th element of the vector $\vec{u}^{(j)}$, which in turn denotes the numerical solution on the $j$th block. The SATs on the right-hand side couple the numerical solution across block interfaces in much the same way as numerical fluxes are used in finite element methods. The scalar $\theta$ is a yet unspecified parameter. The scheme is conservative and if $\theta \geq 0$, it satisfies the discrete energy estimate (see \cite{eriksson2011stable})
\begin{equation} \label{eq:energy_estimate}
\frac{\text{d}}{\text{d}t} \sum_{k=1}^K \| \vec{u}^{(k)} \|_{P_k}^2 = -\theta \sum_{k=1}^K (u_0^{(k)} - u_N^{(k-1)})^2 \leq 0.
\end{equation}
Consequently, it is also stable. Here, for any vector $\vec{v}$, we have defined $\| \vec{v} \|^2_{P_k} = \vec{v}^\top P_k \vec{v}$. Due to the periodicity of the problem, we use the convention $u_N^{(0)} \equiv u_N^{(K)}$. Throughout, we will use $\theta = 1$.

Suppose that the above scheme is used up to some time $t^* > 0$, after which we change from the SBP operator $(P_k,Q_k)$ to some other SBP operator $(P_k^*,Q_k^*)$. Let $\vec{u}^*$ denote the numerical solution after the change. Imposing that $\vec{u}^*(t^*) = \vec{u}(t^*)$, this change constitutes an example of a so-called transmission problem \cite{nordstrom2018well}. A necessary and sufficient condition for retaining the energy estimate \eqref{eq:energy_estimate} through the operator shift is that
\begin{equation} \label{eq:transmission}
    P_k - P_k^* \geq 0.
\end{equation}
Recalling that $P_k$ and $P_k^*$ are diagonal, the following result holds:

\begin{proposition}
Consider two consistent SBP operators $(P_k,Q_k)$ and $(P_k^*,Q_k^*)$. A necessary and sufficient condition for the estimate \eqref{eq:energy_estimate} to be retained after the operator shift is that $P_k = P_k^*$.
\end{proposition}

\begin{proof}
Sufficiency is immediate. For necessity, note that the matrix $P_k - P_k^*$ is diagonal, hence its eigenvalues are given by the diagonal elements. The diagonal elements of both $P_k$ and $P_k^*$ constitute weights in quadrature rules that integrate constants exactly \cite{linders2018order}. Hence,
\begin{equation} \label{eq:summation}
    \vec{1}^\top (P_k - P_k^*) \vec{1} = 0,
\end{equation}
where $\vec{1} = (1,\dots,1)^\top$. Condition \eqref{eq:transmission} dictates that the diagonal elements of $P_k - P_k^*$ are non-negative whereas \eqref{eq:summation} states that their sum is zero. This can only simultaneously hold if $P_k = P_k^*$.
\end{proof}

In subsequent sections we consider adaptive SBP operators that may be inconsistent. Nonetheless, experiments suggest that choosing $P_k = P_k^*$ helps to ensure a robust scheme, hence we will do so throughout.

We remark that an alternative implementation of transmission problems was suggested in \cite{lundquist2020stable}. There, it was shown that stability is retained by imposing $\vec{u}^*(t^*) = 2(I + P_k^{-1} P_k^*)^{-1} \vec{u}(t^*)$. Since the matrices involved are diagonal and positive definite, the inverse exists and is easy to compute. This approach permits $P_k \neq P_k^*$ without stability issues, which unlocks additional degrees of freedom in the optimzation problem described in the next section. However, a detailed study of this approach is beyond the scope of this note.


\section{Stencil optimization} \label{sec:optimization}

Let $\vec{u}$ be a grid function and $\vec{v} \approx \vec{u}_x$ be an accurate approximation of its derivative. We are looking for an SBP operator $D = P^{-1}Q$ that minimizes $D \vec{u} - \vec{v}$. We concentrate on an operator with the following sparsity pattern:
\begin{align*}
    P &= \dx \, \diag{p_0,p_1,p_2,p_3,1,\dots}, \\
    Q &=
    \begin{pmatrix}
        -\frac{1}{2} & q_{0,1} & q_{0,2}  & q_{0,3} &   &   &   &   & \\
        -q_{0,1} & 0 & q_{1,2} & q_{1,3}  &   & &   &   & \\
        -q_{0,2} & -q_{1,2} & 0 & q_{2,3} & a_2 &   &   &   & \\
        -q_{0,3} & -q_{1,3} & -q_{2,3} & 0 & a_1 & a_2 &   &   & \\
        &  & -a_2 & -a_1  & 0 & a_1 & a_2 &  & \\
        & &  &  &  & \ddots &  &  &  \\
    \end{pmatrix}.
\end{align*}
As a reference, we use the conventional SBP operator SBP(4,2), with stencil coefficients
\begin{align*}
    (p_0,p_1,p_2,p_3) &= (17/48,59/48,43/48,49/48), \\
    (a_1,a_2) &= (2/3, -1/12), \\
    (q_{0,1},q_{0,2},q_{1,2},q_{0,3},q_{1,3},q_{2,3}) &= (59/96, -1/12, 59/96, -1/32, 0, 59/96).
\end{align*}
As pointed out in the previous section, it is beneficial from a stability standpoint to fix $P$ for the adaptive SBP operator. Since the target function $\vec{u}$ could be a quadratic polynomial, and since SBP(4,2) is optimal in this case, we will choose the corresponding $P$ also for the adaptive operator. In practice we therefore seek to minimize $Q \vec{u} - P \vec{v}$, where the coefficients of $Q$ constitute the unknowns.

The minimization problem can be rearranged as $A(\vec{u}) \cvec{w} = \vec{b}(\vec{u},\vec{v})$ where the vector $\cvec{w}$ contains the unknowns in $Q$. In our case, $A(\vec{u}) \in \mathbb{R}^{(N+1) \times 14}$ since $Q$ contains $14$ unique coefficients\footnote{Six coefficients per boundary and two in the interior.}. This constitutes a \emph{rank deficient least squares problem}, i.e. there are fewer unknowns than conditions, yet for any $\vec{u}$ some columns in $A(\vec{u})$ are linearly dependent.

Conventional SBP operators are derived by imposing a sequence of such minimization problems where $\vec{u}$ and $\vec{v}$ are substituted by polynomials of increasing degree. For example, SBP(4,2) is uniquely determined by imposing
\[
\begin{aligned}
&\text{first } A(\vec{1}) \cvec{w} = \vec{b}(\vec{1},\vec{0}), \\
&\text{then } A(\vec{x}) \cvec{w} = \vec{b}(\vec{x},\vec{1}), \\
&\text{then } A(\vec{x}^2) \cvec{w} = \vec{b}(\vec{x}^2,2\vec{x}).
\end{aligned}
\]
Each of these conditions can be satisfied exactly, i.e. SBP(4,2) differentiates quadratic polynomials without error.

Suppose that the vector $\vec{u}$ is the numerical solution of a PDE at some point in time and that $\vec{v} \approx \vec{u}_x$ is an approximation of its derivative. The adaptive SBP operator considered here is uniquely defined by solving the sequence of least squares problems
\[
\begin{aligned}
&\text{first } A(\vec{u}) \cvec{w} = \vec{b}(\vec{u},\vec{v}), \\
&\text{then } A(\vec{1}) \cvec{w} = \vec{b}(\vec{1},\vec{0}), \\
&\text{then } A(\vec{x}) \cvec{w} = \vec{b}(\vec{x},\vec{1}).
\end{aligned}
\]
Thus, we first attempt to choose coefficients that accurately differentiate $\vec{u}$, then use the remaining degrees of freedom to impose conventional accuracy conditions. An important difference to SBP(4,2) is that none of these conditions can be satisfied exactly; they are all imposed in a least squares sense.


\section{Numerical experiments} \label{sec:experiments}

We explore the efficacy of the adaptive SBP operator by solving the periodic advection problem \eqref{eq:advection} on $K$ elements with $N+1$ grid points each. As the exact solution we choose $u(x,t) = \sin{(2 \pi (x-t))} + \cos{(4 \pi (x-t))}/2$.

To obtain the vector $\vec{v}$ needed to optimize the SBP operator, we use a block norm SBP operator for which $Q$ has identical sparsity pattern to SBP(4,2), and $P$ has a $4 \times 4$ block near each boundary. Its boundary accuracy is $\ordo{3}$.

To obtain any accuracy at all with the adaptive SBP operators it is necessary to re-optimize sufficiently often. This is a consequence of the adaptive operator being an accurate derivative approximation for $\vec{u}$ only at a particular time. As $\vec{u}$ evolves in time, the quality of any given operator deteriorates. Here we make an empirical choice to re-optimize at time intervals $\Delta \tau = 1/(2 K (N+1))$.

Time integration is performed using Matlab's ode45 routine (Dormand-Prince) with absolute and relative tolerances set to $10^{-10}$. It is thus expected that the time integration has negligible impact on the overall error.

\begin{figure}
     \centering
     \begin{subfigure}[b]{0.49\textwidth}
         \centering
         \includegraphics[width=\textwidth]{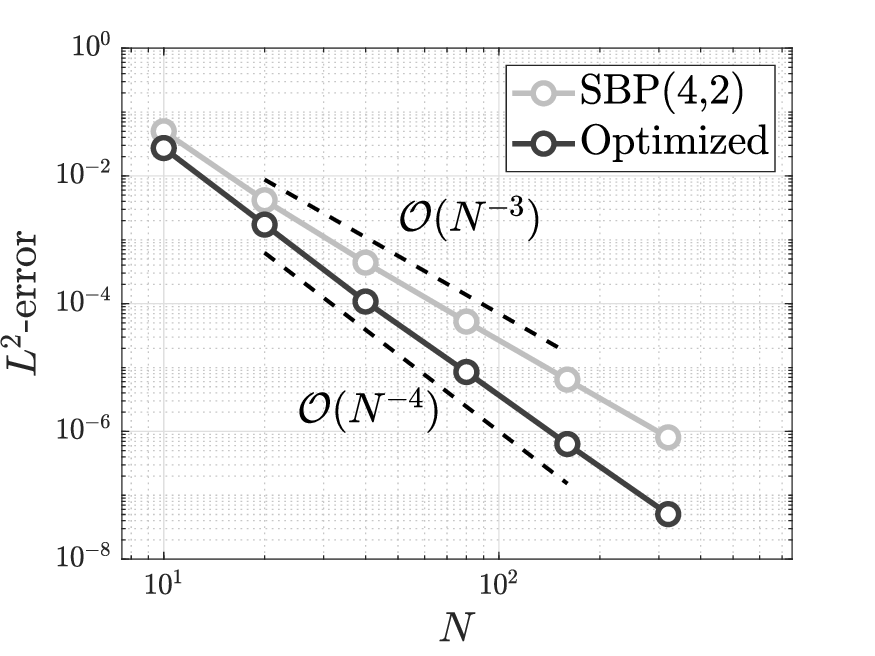}
         \caption{Grid convergence}
         \label{fig:N_conv_4}
     \end{subfigure}
     \hfill
     \begin{subfigure}[b]{0.49\textwidth}
         \centering
         \includegraphics[width=\textwidth]{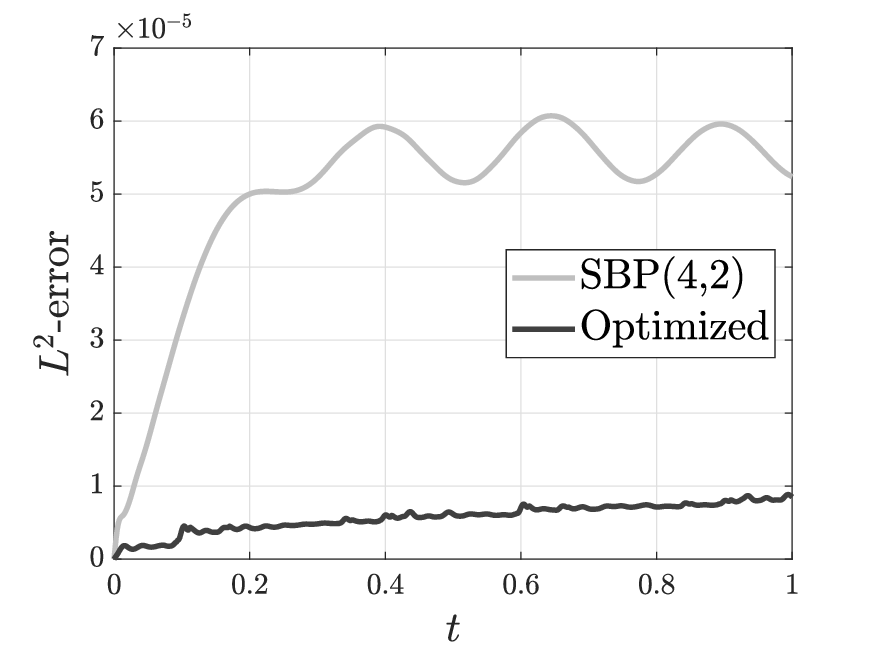}
         \caption{Time error}
         \label{fig:t_error_4}
     \end{subfigure}
    \caption{$L^2$-errors with $K=4$.}
    \label{fig:K_4}
\end{figure}

We set $K=4$ so that four SBP operators make up the spatial discretization. These are independently optimized, each using a quarter of the information contained in $\vec{u}$. Fig. \ref{fig:N_conv_4} shows the convergence of the $L^2$-errors of SBP(4,2) and the adaptive method when $N$ is varied. Fig. \ref{fig:t_error_4} shows the time evolution of the errors when $N=80$. On all grids the errors are visibly smaller with the adaptive SBP operators. On the finest grids, the error is reduced by more than an order of magnitude. The observed convergence rate is $\ordo{3}$ for SBP(4,2) as expected\footnote{Here, $\dx$ is a constant multiple of $N^{-1}$ shown in Fig. \ref{fig:N_conv_4}}. Remarkably, for the adaptive stencil it is $\ordo{4}$. Thus, the adaptive stencil overcomes the order limit for non-adaptive SBP operators. This happens despite the optimized SBP operator formally not having any order at all. It is not even consistent.

Repeating the experiment with $K=1$ elements and four times as many grid points per element yields very similar results to those in Fig. \ref{fig:K_4} (not shown). Since a single SBP operator now makes up the discretization, there are fewer stencil coefficients to optimize. On the other hand, that optimization has access to all the data in $\vec{u}$.


\section{Generalizations} \label{sec:generalization}

We have considered the special case of an adaptive SBP operator whose sparsity pattern matches SBP(4,2). Several extensions are required to turn the approach into a general method:

\begin{itemize}
    \item Extending the adaptive method to SBP operators using central differences of bandwidth $2n+1$ and boundary blocks of size $r \times r$ for $n>2$ and $r>4$ unlocks additional free parameters. A longer sequence of least squares problems must be solved, approximately imposing exactness for higher degree polynomials. However, a complication arises where, numerically, the over-determined linear systems have some singular values near machine epsilon. The optimization is highly sensitive to these spurious singular values, which must be identified and deflated. Already with $(n=3,r=6)$, standard least squares solvers fail for the system $A(\vec{x}^3) \cvec{w} = \vec{b}(\vec{x}^3,3\vec{x}^2)$.

    \item Additional desirable properties such as a small spectral radius can be imposed using appropriate constraints. 

    \item The selection of re-optimization times should be automized, which requires the development of a sensor.

    \item Solving a sequence of least squares problems comes with a non-negligible expense. We have made no attempt at optimizing the procedure, however this will be necessary for the benefits to outweigh the cost of adaptivity.
\end{itemize}


\section{Conclusions} \label{sec:conclusions}

We have experimentally illustrated the improved accuracy of a stencil-adaptive SBP operator with a sparsity pattern matching SBP(4,2). The most notable feature is its superconvergence. While we have not provided a theoretical basis for this behavior, we remark that, to the best of our knowledge, this is the first observation of such rapid convergence for SBP finite difference methods with diagonal $P$.


\bibliographystyle{model1-num-names}
\bibliography{refs}

\end{document}